\documentclass[final,3p,times]{elsarticle}

\usepackage{lineno,hyperref}
\modulolinenumbers[5]                

\usepackage{amsfonts,amssymb,amsmath,amsthm}
\usepackage{graphicx}
\usepackage{color}

\newcommand{\hX}{{\widehat{X}}}

\newcommand{\C}{\mathcal{C}}

\renewcommand{\S}{\mathcal{S}}
\newcommand{\U}{\mathcal{U}}

\newtheorem{theorem}{Theorem}

\newtheorem{proposition}[theorem]{Proposition}

\theoremstyle{remark}
\newtheorem{remark}{Remark}

\usepackage{amsfonts}

\begin{document}

\begin{frontmatter}

\title{Local stability in a transient Markov chain}

\author[TUe]{Ivo Adan}
\author[HWUNSU]{Sergey Foss}
\author[HWU]{Seva Shneer}
\author[UH]{Gideon Weiss}

\address[TUe]{Eindhoven University of Technology, Eindhoven, the Netherlands}
\address[HWUNSU]{Heriot-Watt University, Edinburgh, UK and MCA, Novosibirsk State University and Sobolev Institute of Mathematics, Novosibirsk, Russia}
\address[HWU]{Heriot-Watt University, Edinburgh, UK}
\address[UH]{The University of Haifa, Haifa, Israel}


\begin{abstract}
We prove two propositions with conditions that a system, which is described by a transient Markov chain, will display local stability.  Examples of such systems include partly overloaded Jackson networks, partly overloaded polling systems, and overloaded multi-server queues with skill based service, under first come first served policy.
\end{abstract}


\end{frontmatter}



\section{Introduction}
\label{sec.introduction}
Many complex stochastic systems can be described by an irreducible Markov chains on a countable state space.  It is often the case that the state of this Markov chain is composed of several components, where each component describes the ``local'' state of part of the system.  While the dynamics of the system, given by the transition mechanism of the Markov chain, are influenced by the state of the entire system, it is often the case that the dependence of the local transitions is only weakly coupled with the rest of the system.

Essential to the study of  Markov chains is the question of stability:  Is the chain ergodic, in which case it has a stationary distribution from which its long time average behavior can be obtained, or is it transient, in which case it may be studied through fluid approximations.  These two modes of behavior are totally different.

In complex systems one may however be faced by an intermediate sort of behavior.  While the system as a whole is transient, and so there is no stationary distribution for the entire system, some components of the system, when regarded locally, display stable behavior and seem to approach a stationary distribution when regarded on their own, at least for most of the time.  

We formulate this type of situation, and prove that under the proper conditions one can indeed talk about local stability of such transient systems.  This is done in Propositions \ref{thm.lemma} and \ref{l2}, in Section \ref{sec.lemma}.  Before that, in Section \ref{sec.examples}, we present some examples.

We point out also that the setting here is similar to the situation where some components of a Markov chain are positive recurrent Markov chains on their own, and the limiting behaviour of other components, as well as that of the entire process, may be computed through averaging over the stationary distribution of the stable component(s). Such examples have been considered in, among other work, \cite{ethier2009markov,szpankowski1988stability,borst2008stability,foss2013stability}.


\section{Examples of local stability}
\label{sec.examples}
  
\subsection{Jackson networks}
\label{sec.jackson}

This example is analyzed by Goodman and Massey \cite{goodman-massey:84}.
In a Jackson network, single server nodes $i=1,\ldots,I$ have exogenous arrival rate $\nu_i$, service rate $\mu_i$, and routing probabilities $P_{i,j}$ for a customer that completes service at node $i$ to go next to node $j$, with the matrix $P$ sub-stochastic with spectral radius $< 1$.  The traffic equations for a stable Jackson network are:
$
\lambda = \nu + P' \lambda,
$
with $\lambda_i$ the stationary rate of inflow and outflow of customers from node $i$.
Necessary and sufficient for ergodicity is $\lambda_i<\mu_i$, i.e. $\rho_i = \lambda_i/\mu_i<1$, for $ i=1,\ldots,I$.
If the ergodicity condition does not hold for all nodes, then the system is transient.  The traffic equations are  now modified to
$
\lambda = \nu + P' (\lambda \wedge \mu),
$
which has a unique solution, reached by solving a linear complementarity problem.  It divides the nodes into two sets:  $I_0=\{i:\lambda_i < \mu_i\}$ and $I_1=\{i:\lambda_i \ge \mu_i\}$. 
An intuitive idea (made rigorous in \cite{goodman-massey:84}) here is that nodes in $I_1$ diverge (are unstable), 
while all the nodes in $I_0$ act like a Jackson network with augmented input, where for each $i \in I_0$, in addition to exogenous input rate $\nu_i$,  there is input from each node $j\in I_1$, at rate $\mu_j P_{j,i}$.

The difficulty in verifying the above is that there is no steady state distribution for the whole system, since it is transient, and there are no equilibrium equations for the subsystem of node $I_0$, since they do not form a Markov chain.  Goodman and Massey prove the statement by considering two ergodic Markov chains for $I_0$ which provide stochastic upper and lower bound for $Q_i(t),\,i\in I_0$, and one of which has the exact limiting stationary distribution, while  the other has a parameter $\epsilon$,  and its stationary 
distribution converges to the same limiting one as $\epsilon \to 0$.


 \subsection{Multi-server queues with skill based service under FCFS policy}
 
 This example is analyzed by Adan and Weiss \cite{adan-weiss:12} (see also \cite{adan-weiss:11}) .  In a skill based service queue  there are  customers of types $ \C=\{c_1,\ldots,c_I\}$ and servers $\S=\{s_1,\dots,s_j\}$, and a bipartite compatibility graph between $\S$ and $\C$, with an arc $(s_j,c_i)$  if server $s_j$ can serve customers of type $c_i$.  Assume customer arrivals are Poisson at rates $\lambda_{c_i}$, and service is exponential, with rates $\mu_{s_j}$, let $\lambda=\sum_\C \lambda_{c_i}$, $\mu=\sum_\S \mu_{s_j}$, 
 $\alpha_{c_i}=\lambda_{c_i}/\lambda$, $\beta_{s_j} = \mu_{s_j}/\mu$.  
Denote $\S(c_i)$ the servers of $c_i$, $\C(s_j)$ the customers of $s_j$, and for  $C\subseteq \C$, $S\subseteq \S$,  let $\S(C)=\bigcup_{c_i\in C} \S(c_i)$, $\C(S)=\bigcup_{s_j\in S} \C(s_j)$, and let also $\U(S)=
\overline{\C(\overline{S})}$ be  customer types which can only be served by servers in $S$.  
Let $\alpha_C=\sum_{c_i\in C} \alpha_{c_i}$ and $\beta_S=\sum_{s_j\in S} \beta_{s_j}$, with analogous notation for $\lambda_C,\,\mu_S$.

Service discipline is first come first served (FCFS) assign longest idle server (ALIS),  i.e. server $s_j$, when free, will take the longest waiting compatible customer, and arriving customer of type $c_i$ is assigned to longest idle compatible server.


The authors of \cite{adan-weiss:12} construct a Markov chain $X(t)$ describing the state of the system and show that this Markov chain is ergodic if and only if for  every non-empty subset of customer types $C$, and of servers $S$, the three equivalent sets of conditions hold:
\begin{equation}
\label{eqn.cond1}
\lambda_C <  \mu_{\S(C)},  \qquad  \mu_S < \lambda_{\C(S)}, \qquad \mu_S > \lambda_{\U(S)}.
\end{equation}


If condition (\ref{eqn.cond1}) fails, then $X(t)$ is transient, and  queues of some types of customers will grow to infinity.  However, local stability, in the sense that some components of $X(t)$ converge to an intuitive stationary distribution, may still hold (see \cite{adan-weiss:12} for details).

\subsection{Mesh network governed by a CSMA/CA protocol}

A further example was considered in \cite{Shneer} and \cite{ShneerVen2015} (see also references therein) and concerns a wireless network consisting of a number of nodes placed on a line. Messages enter the system at the left-most node and need to be relayed from node to node until they are transmitted by the right-most node, at which point they leave the network. A complication is that neighbouring nodes cannot transmit simultaneously, due to interference, and the so-called CSMA/CA algorithm is used to regulate which nodes transmit at which time. It is well known (clear intuitively and shown in experiments) that the queues of some nodes grow to infinity, yet other queues remain stable. A rigorous analysis found in \cite{ShneerVen2015} leads to the characterisation of the end-to-end throughput of an example of such a network.

\subsection{Polling Systems}
\label{sec.polling}
Foss, Chernova and Kovalevskii \cite{foss-etal:96} consider a polling system with one or several servers, and a number of stations.  Input consists of stationary ergodic streams of customers. Servers follow i.i.d. cyclic routes through all the stations,  which are independent of the arrivals and of the queues. The service policy has a number $f_k^j(x,D_k^j)$ of customers served on the $j$th visit of the server to station $k$, if there are $x$ customers at the station, where $\{D_k^j\}_{k,j}$ form an i.i.d. sequence, and service is monotone in the sense that $f_k(x,D) \le f_k(x+1,D) \le f_k(x,D)+1$. In these systems it is possible that some of the queues are unstable while 
the remaining queues display local  stability.


\section{Local stability of transient Markov chains}

For the sake of clarity, we formulate our results in a particular
case of a discrete-time Markov chain on a countable state space.
Out results continue to hold in much more general settings, see Remark 1 below for more details. 

\label{sec.lemma}
\begin{proposition} 
\label{thm.lemma}
Let $X(n)=(X_1(n),X_2(n))$ be a Markov chain on a countable state space with \\$X_1(n) =(X_{11}(n),\ldots,X_{1,k}(n))\in
\mathbb{Z}^k_+$ and $X_2(n)=(X_{21}(n),\ldots,X_{2m}(n))\in \mathbb{Z}^m_+$. Assume the following:
\begin{enumerate}
\item $\lim_{n \to \infty} X_{2i}(n) = \infty$ almost surely, for all $i=1,\ldots,m$ and for any initial condition
$(X_1(0),X_2(0))$.
\item ${\mathbf P}(X_1(n+1)=j |
    X_1(n)=i,\;X_2(n)=l) = P_{i,j}$, for all values of $l=(l_1,\ldots,l_m)$ with all strictly positive coordinates, where $P_{i,j}$ are
    transition probabilities of an ergodic Markov chain with the unique stationary
    distribution $\pi =\{\pi_j\}$.
\end{enumerate}
Then for all initial  $i_0,j_0$:
\[
\sup_j  \Big{|}{\mathbf P}\big{(}X_1(n)=j \,|\, X_1(0)=i_0,X_2(0)=j_0\big{)} - \pi_j \Big{|} \to 0, 
\mbox{ as } n \to  \infty,
\]
i.e. $X_1(n)$ converges in distribution to $\pi$ in the total variation norm.
\end{proposition}
\begin{proof} To simplify the notation, we provide a proof for $k=m=1$. 
Denote the transition probabilities of $X(t)$ as follows:
\[
{\mathbf P}(X_1(n+1)=k,X_2(n+1)=l\,|\, X_1(n)=i, X_2(n)=j) = P_{(i,j),(k,l)}.
\]
By Condition 2 we can write, whenever $j>0$:
\begin{align*}
& {\mathbf P}(X_2(n+1)=l\,|\, X_1(n+1)=k,X_1(n)=i, X_2(n)=j) \\ & =\frac{{\mathbf P}(X_1(n+1)=k,X_2(n+1)=l\,|\, X_1(n)=i, X_2(n)=j)} {{\mathbf P}(X_1(n+1)=k\,|\, X_1(n)=i,
X_2(n)=j)} = \frac{P_{(i,j),(k,l)}}{P_{i,k}}.
\end{align*}

Fix initial values $X_1(0)=i_0,X_2(0)=j_0$, and choose arbitrary $\epsilon >0$.
Choose large enough $n_0$ (to be specified).  Starting at time $n_0$ with the values
$X_1(n_0),X_2(n_0)$, we now construct a new Markov chain, $(Y(n),\hX_1(n),\hX_2(n))$
for $n=n_0,n_0+1,\ldots$. We initialize them at time  $n_0$ as
$Y(n_0) = \hX_1(n_0)=X_1(n_0)$ and $\hX_2(n_0) =  X_2(n_0)$.

Starting from these values at $n_0$ the following transitions are made, from $n$ to
$n+1$ for $n\ge n_0$:
\begin{eqnarray*}
&& {\mathbf P}(Y(n+1) = k\,|\,Y(n)=i) = P_{i,k}, \\ && \hX_1(n+1) =  \begin{cases}
Y(n+1) & \mbox{ if } \hX_2(n) > 0 \\  {x} & \mbox{ if } \hX_2(n) = 0,
\end{cases} \\
\end{eqnarray*}
with an arbitrary $x$. The value of $\hX_2(n+1)$ is generated as follows:
\begin{eqnarray*}
&&{\mathbf P}(\hX_2(n+1) = 0\,|\, Y(n)=k,\hX_1(n)=i,\hX_2(n)=0)=1, \\ && 
{\mathbf P}(\hX_2(n+1)=l\,|\,
\hX_1(n+1)=k,\hX_1(n)=i, \hX_2(n)=j, j>0) = \frac{P_{(i,j),(k,l)}}{P_{i,k}}.
\end{eqnarray*}

Observe that $Y(n)$ on its own is a Markov chain, with transition probabilities
$P_{i,k}$.  Also, $\hX_1(n),\hX_2(n)$ on its own is distributed for $n\ge n_0$ exactly
like $X_1(n),X_2(n)$, for as long as $\hX_2(n-1)> 0$.  Once $\hX_2(n)=0$, it will stay
as $0$ for all times $m\ge n$, and $\hX_1(m)=x$ for all $m>n$.  Finally, note that
$\hX_1(n)=Y(n)$ for as long as $\hX_2(n-1) > 0$.

The following holds:
\[
{\mathbf P}(X_1(n) = j \mbox{ and } X_2(m)>0 \mbox{ for all }n_0 \le m \le n-1) =
{\mathbf P}(\hX_1 (n) = j) = {\mathbf P}(\hX_1 (n) = Y(n) = j)
\]

We now look at the total variation distance between the distribution of $X_1(n)$
(having started from the fixed $i_0$ at time 0), and the distribution $\pi$.
\begin{align*}
& | {\mathbf P}(X_1(n) = j) - \pi_j | = | {\mathbf P}(X_1(n) = j \mbox{ and } X_2(m)>0 
\mbox{ for all }n_0
\le m \le n-1) \\ & + {\mathbf P}(X_1(n) = j 
\mbox{ and }
X_2(m)=0\mbox{ for some } n_0 \le m \le n-1) -  \pi_j |\\ & = |
{\mathbf P}(\hX_1(n)=Y(n)=j) + {\mathbf P}(X_1(n) = j \mbox{ and } X_2(m)=0\mbox{ for 
some } n_0 \le m \le
n-1) - \pi_j |\\ &  = | {\mathbf P}(Y(n)=j) - \pi_j - {\mathbf P}(Y(n)=j, \hX_1(n)\ne Y(n))
 +{\mathbf P}(X_1(n) = j \mbox{ and } X_2(m)=0\mbox{ for some } n_0 \le m \le n-1) | \\
& \le | {\mathbf P}(Y(n)=j) - \pi_j| + {\mathbf P}(Y(n)=j, \hX_1(n)\ne Y(n))
 + {\mathbf P}(X_1(n) = j \mbox{ and } X_2(m)=0\mbox{ for some } n_0 \le m \le n-1)  \\
&  \le | {\mathbf P}(Y(n)=j) - \pi_j| +2 {\mathbf P}( X_2(m)=0\mbox{ for some } 
n_0 \le m \le n-1)
 \\&  \le | {\mathbf P}(Y(n)=j) - \pi_j| +2 {\mathbf P}( X_2(m)=0\mbox{ for some } n_0 \le m <
\infty).
\end{align*}

We now make use of Condition 1 to choose $n_0$. Since $\lim_{n \to \infty} X_2(n) =
\infty$ almost surely, we have for every $j$ a random time (not a stopping time)
$\nu(j)$ such that
$\nu(j) = \sup \{n: X_2(n) < j \}$ and ${\mathbf P}(\nu(j) < \infty) = 1$.  We can now choose $n_0$ large enough so that
${\mathbf P}( \nu(1) \ge n_0) < \epsilon/2$.
Hence, with probability exceeding $1-\epsilon/2$, $X_2(n)>0$ for all $n \ge n_0$, and
we have
$| {\mathbf P}(X_1(n) = j) - \pi_j | < | {\mathbf P}(Y(n)=j) - \pi_j| + \epsilon$.

We therefore have:
\begin{equation*}
\limsup_{n\to\infty}  \left( \sup_j | {\mathbf P}(X_1(n) = j) - \pi_j | \right) 
\le \limsup_{n\to\infty}  \left( \sup_j | {\mathbf P}(Y(n)=j) - \pi_j| \right) + \epsilon,
\end{equation*}
but
$\limsup_{n\to\infty}  \left( \sup_j | {\mathbf P}(Y(n)=j) - \pi_j| \right) = \lim_{n\to\infty}
\left( \sup_j | {\mathbf P}(Y(n)=j) - \pi_j| \right) = 0$
and we have shown that
$$
 \limsup_{n\to\infty}  \left( \sup_j | {\mathbf P}(X_1(n) = j) - \pi_j | \right)  \le \epsilon
$$
for an arbitrary $\epsilon >0$.  This completes the proof.
\end{proof}

\begin{remark}  
As stated above, the same scheme works in a more general setting, of a general measurable
state space Markov process $\mathcal{X}(t)$,  in discrete  time,  where
assumption (1) is replaced by the assumption that some test function $L(
\mathcal{X}(t)) \to \infty$ almost surely as $t\to\infty$, and assumption (2) says that
conditional on $L(\mathcal{X}(t))>0$, a process which is a function of $\mathcal{X}(t)$, say $M(\mathcal{X}(t))$, satisfies that
$M(\mathcal{X}(t))\,|\,L(\mathcal{X}(t))=l$ is ergodic, and independent of the value
$l>0$.  The conclusion then is that as $t\to\infty$,  the distribution of $M(\mathcal{X}(t))$ converges
to the invariant distribution of $M(\mathcal{X}(t))\,|\,L(\mathcal{X}(t))=l>0$.  A similar conclusion holds in continuous time.
\end{remark}

\begin{remark}

If one thinks of $X_2(n)$ as being an autonomous Markov chain, then Condition 1 of the proposition above means that $X_2(n)$ is transient. A natural question is what happens if we replace this condition with the requirement that $X_2(n)$ is null-recurrent. The results of \cite{goodman-massey:84} hold in the case that some of the nodes in a Jackson network have service rates which are exactly equal to the arrival rates. The question arises whether this could be proven for a general two-component  state Markov chain, i.e. whether the proposition above is valid if one replaces Condition 1 with the convergence of $X_2(n)$ to infinity in distribution.

The proofs in \cite{goodman-massey:84} are based on the monotonicity of the network, and the convergence of the distribution of $X_1(n)$ to the natural limiting one will indeed hold for all specific models exhibiting such monotonicity. A relatively general statement (Proposition 2) and its proof are given below. However, such a convergence does not hold in general, and here is an example. 

{\bf Example}.
Let $X_2(n)$ be a simple random walk on $\mathbb{Z}_+$ reflected at $0$:
$$
X_2(n) = \max\{0, X_2(n-1) + \xi(n)\},
$$
where $\{\xi(i)\}$ are i.i.d with ${\mathbf P}(\xi(1)=1) = {\mathbf P}(\xi(1) = -1) = 1/2$ and $X_2(0)=0$. 

Define the first moment this random walk returns to $0$ as
\begin{equation} \label{eq:hitting_x2}
t(1) = \inf\{n \ge 1: X_2(n) = 0\}.
\end{equation}
It is known that $t(1) < \infty$ a.s. but $\mathbf{E} t(1) = \infty$. In fact ${\mathbf P}(t(1) =k) \sim C k^{-3/2}$ as $k \to \infty$. Assume now that
$
X_1(n) = \max\{0, X_1(n-1) + \eta(n)\},
$
where
\begin{eqnarray*}
\eta (n) =
\begin{cases}
-1, \quad \mbox{if} \quad  X_2(n-1) > 0,\\
\psi(n), \quad \mbox{if} \quad X_2(n-1) = 0
\end{cases}
\end{eqnarray*}
for an i.i.d sequence $\psi(n)$ that does not depend on the dynamics of $\{X_2(n)\}$ and is
such that ${\mathbf P}(\psi(n) > k) \sim k^{-\alpha}$ as $k\to\infty$, where $0< \alpha < 1/2$.  
The distributions of $t(i)$ and $\psi(i)$ are regularly varying and, therefore,
{\it long-tailed}. Then (see e.g. Chapter 2 of \cite{FKZ}), 
\begin{equation}\label{equi}
{\mathbf P} (\psi(1)-t(1) >x) \sim {\mathbf P} (\psi(1)>x)
\quad \mbox{and}
\quad
{\mathbf P} (t(1)-\psi(1)>x)\sim {\mathbf P} (t(1)>x) \quad \text{as} \quad x\to \infty.
\end{equation}

The evolution of $X_1(n)$ is rather simple: the value is decremented by $1$ at each step while $X_2(n)$ is positive and jumps up by a random variable with a distribution of $\psi(1)$ when $X_2(n)$ hits zero.
It is clear that, conditioned on $X_2(n)$ staying always positive, $X_1(n)$ converges to $0$, regardless of its starting point. 
We are going to show now that, unconditionally,  
$X_1(n)$
converges to infinity a.s. For that, it is enough
to show that random variable
$
\tau = \inf\{n \ge 1: X_1(n) = 0\}
$
is improper, i.e. 
${\mathbf P}(\tau = \infty) > 0$. 
Along with (\ref{eq:hitting_x2}), define
$
t(k+1) = \inf\{n > t(1)+\ldots +t(k): X_2(n) = 0\} - \sum_{j=1}^k t(j)$ for $k \ge 1$.
Define now a random walk
$
S_n = \sum_{i=1}^n (\psi(i) - t(i)).
$
Assume that $X_1(0)=X_2(0)=0$. Then  $\zeta(i) := \psi(i)-t(i)$ are i.i.d. random variables, with ${\mathbf E} \zeta^+ = {\mathbf  E} \zeta^-=\infty$
where $\zeta^+ = \max (\zeta(1),0)$ is the positive part and $\zeta^- = \max (-\zeta(1),0)$ the negative part of random variable $\zeta(1)$. It is clear that
$$
{\mathbf P}(\tau = \infty) = {\mathbf P}(X_1(n) > 0 \quad \mbox{for all} \quad n> 0) = {\mathbf P}(S_n > 0 \quad \mbox{for all} \quad n> 0)
$$
and it is sufficient to show that $S_n \to \infty$ a.s.  
For this we apply Theorem 2 of \cite{erickson:73}. In order to be consistent with the notation of the paper, we let 
$
m_+(x) = {\mathbf E} \min (\zeta^+,x)
$
, for $x>0$. By Theorem 2 of \cite{erickson:73}, $S_n/n\to\infty$ a.s. (and then $S_n\to\infty$ a.s.) 
if and only if 
$$
J_- = {\mathbf E} \frac{\zeta^-}{m_+(\zeta^-)} < \infty.
$$
By \eqref{equi}, $m_+(x) \sim {\mathbf E} \min (\psi(1), x) \sim x^{1-\alpha}/(1-\alpha )$ as $x\to\infty$ and, further,
$
J_- \approx (1-\alpha ){\mathbf E} \frac{t(1)}{t(1)^{1-\alpha}} = (1-\alpha ) {\mathbf E} t(1)^{\alpha} < \infty
$
(recall that ${\mathbf P}(t(1) =k) \sim C k^{-3/2}$ as $k \to \infty$).
\end{remark}

\begin{remark}
It is clear that the result of the proposition above holds if one replaces its Condition 2 by the the requirement that
${\mathbf P}(X_1(n+1)=j |
    X_1(n)=i,\;X_2(n)=l) = P_{i,j}$, for all values of $l=(l_1,\ldots,l_m)$ with all $l_i>N$,
for any fixed $N$. An interesting question is to consider the case when
$
{\mathbf P}(X_1(n+1)=j |
    X_1(n)=i,\;X_2(n)=l) \to P_{i,j}
$
when $l \to \infty$, i.e. making the dynamics of $X_1(n)$ asymptotically independent of the position of $X_2(n)$ rather than simply independent of it.
This question requires further efforts, and we plan to pursue this research direction.

\end{remark}

\begin{proposition}\label{l2}
Assume again that $X(n)=(X_1(n),X_2(n))$ is a Markov chain taking values in $\mathbb{Z}^{k+m}_+$. Suppose that Condition 2 of Proposition 1 continues to hold, while Condition 1
is replaced by\\
$\widetilde{1}$. $X_2(n)\to\infty$ in probability (again, coordinate-wise), given $X_1(0)=0$ and $X_2(0)=0$.\\
Further, assume that \\
3. Markov chain $X(n)$ is {\it monotone}, in the following sense. For $y\in \mathbb{Z}^{k+m}_+$, let $C_y = \{ z\in \mathbb{Z}^{k+m}_+ \  : \  z\ge y\}$
where $\ge$ is the standard partial ordering in $\mathbb{Z}^{k+m}$. Then {\it monotonicity} means that
$$
{\mathbf P} (X(1)\in C_y \ | X(0)=x) \ge {\mathbf P} (X(1)\in C_y \ | \ X(0)=\widehat{x}) \quad \text{for all} \quad x 
\ge \widehat{x}, \quad y \in \mathbb{Z}^{k+m}_+.
$$


Then the conclusion of Proposition 1 holds again.
\end{proposition}

\begin{proof}
Again, we consider the case $k=m=1$ only.
In what follows, all equalities and inequalities hold a.s.

It is known \cite{BF} that a Markov chain may be represented as a stochastic recursion
$$
X(n+1)= f(X(n), U_n) \equiv (f_1(X(n), U_n), f_2(X(n),U_n)) = (X_1(n+1),X_2(n+1))
$$
where $U_n$ are i.i.d. random variables having the uniform-$(0,1)$ distribution. Also, monotonicity if $X(n)$
implies that $f,f_1,f_2$ may be chosen monotone in their first argument. 
We assume $U_n$ to be given for all $-\infty < n < \infty$.

By assumptions of the proposition, there exists a stationary Markov chain $\widehat{X}_1 (n), -\infty < n <\infty$ with distribution $\pi = \{\pi_j\}$ that is measurable
with respect to $\{U_n\}$ and satisfies recursion
$
\widehat{X}_1(n+1) = f_1((\widehat{X}_1(n),1), U_n), \quad -\infty <n <\infty.
$
Also, for any $l=0,1,2,\ldots$, one may introduce a Markov chain
$\widehat{X}_1^{(-l)}(n), n=-l,-l+1,\ldots$ that starts from $\widehat{X}_1^{(-l)}(-l)=0$ at time $-l$ and satisfies recursion
$\widehat{X}_1^{(-l)}(n+1) = f_1((\widehat{X}_1^{(-l)}(n),1),U_n)$. By the monotonicity,
$$
0 \le \widehat{X}_1^{(-1)}(0) \le \widehat{X}_1^{(-2)}(0) \le \ldots \le \widehat{X}_1(0)
$$
and, moreover, there is an a.s. finite time $\nu$ (called {\it backward coupling time}, see \cite{BF}) such that
$$
\widehat{X}_1^{(-l)}(0) = \widehat{X}_1(0), \quad \mbox{for all} \quad l\ge \nu.
$$
The latter follows from the uniqueness of the stationary distribution and from the fact that
monotone convergence on the lattice must be with coupling (see again \cite{BF} for supportive
arguments). 

Further, for any $l=1,2,\ldots$, consider an
auxiliary Markov chain ${X}^{(-l)}(n)=({X}_1^{(-l)}(n),{X}_2^{(-l)}(n))$, $n=-l,-l+1,\ldots$ that starts from
${X}^{(-l)}(-l)=(0,0)$ at time $-l$. Then, by the monotonicity, for any time $n\le 0$ and for any $-l_1\le -l_2\le n$,
$$
{X}^{(-l_1)}(n) \ge {X}^{(-l_2)}(n)
$$
and also
$
\widehat{X}_1(n)\ge {X}_1^{(-l_1)}(n)\ge {X}_1^{(-l_2)}(n).
$
Then
\begin{equation}\label{monot2}
X^{(-l)}_2(n) \uparrow \infty \quad \mbox{a.s. as} \quad l\to\infty.
\end{equation}
for any fixed $n\le 0$. The sample-path monotonicity  in \eqref{monot2} follows from the fact that the Markov chains
start from the minimal state and from the induction arguments. Indeed, $X^{(-l-1)}(-l) \ge (0,0)= X^{(-l)}(-l)$ a.s
and then, for any $-l \le n$,
if $X^{(-l-1)}(n) \ge X^{(-l)}(n)$ a.s, then 
$
X^{(-l-1)}(n+1) = f( X^{(-l-1)}(n),U_n)\ge f(X^{(-l)}(n),U_n) = X^{(-l)}(n+1) \quad \mbox{a.s.}
$ 
Then the a.s convergence to infinity in \eqref{monot2} 
follows from the monotonicity and the convergence to infinity in probability (since $X_2(n+l)$ and $X^{(-l)}_2(n)$
have the same distribution, for any $n$ and any $l$). Take any small $\varepsilon >0$ and choose $N>>1$ such that ${\mathbf P} (\nu > N)\le \varepsilon/2$
and then $L>N$ such that 
$
{\mathbf P} (X^{(-L)}_2(n)\ge 1, \ \ \mbox{for all} \ \  -N\le n \le -1)\ge 1- \varepsilon/2.
$
Then, on the event
$
A:= \{\nu \le N\} \cap \cap_{n=-N}^{-1} \{X^{(-L)}_2(n)\ge 1\}
$
of probability at least $1-\varepsilon$, we have that
$
X^{(-L)}_1(n)=\widehat{X}^{(-L)}_1(n), \quad \mbox{for all} \quad -N\le n \le 0
$
and, in particular,
$
X^{(-L)}_1(0)=\widehat{X}_1^{(-L)}(0)=\widehat{X}_1(0).
$
Therefore, by the monotonicity, on the event $A$ we have
$$
X^{(-l)}_1(0)=\widehat{X}_1(0) \quad \text{for all} \quad l\ge L.
$$

Now start Markov chain $X(n)$ from $X(0)=(0,0)$ at time $0$.
Then, for any $l\ge L$ and for any set $B$,
$$
|{\mathbf P} (X_1(l)\in B) -\pi (B)| =
|{\mathbf P} (X^{(-l)}_1(0)\in B)-\pi (B)|\le \varepsilon.
$$
Since $\varepsilon >0$ is arbitrary, we may conclude that
the distribution of $X_1(n)$ converges to $\pi$ in the total variation 
if the initial value is $(0,0)$.

If the Markov chains starts now from another initial state, say $X(0)=(i,j)$, then 
it may be squeezed between $X_1(n)$ that starts
from $0$ and, say, sequence $\widetilde{X}_1(n+1)= f_1(\widetilde{X}_1(n),1),U_n)$
that starts from $\widetilde{X}_1(0)=i$. Since both boundary sequences couple with the
stationary sequence $\widehat{X}_1(n)$, the result follows.
\end{proof}

\begin{remark}
In Condition 3 of Proposition 2, we require monotonicity in {\it both} components, $X_1$ and $X_2$.
Because of that, convergence to infinity in \eqref{monot2} is monotone too (and, therefore,
it occurs almost surely). Condition 3 is satisfied for Jackson networks. However, this condition was taken just to make the proof simpler.
One may weaken the condition, by assuming monotonicity in the first component only.
Then convergence in \eqref{monot2} holds in probability, but the statement of the proposition
continues to hold -- one needs a slightly more detailed analysis.
\end{remark}

\bibliography{local_stab_bib}

\end{document}